\documentclass[12pt]{article}
\usepackage{amsfonts}
\usepackage{theorem}
\newcommand{\beq}{\begin{equation}}\newcommand{\eeq}{\end{equation}}
\newcommand{\ba}{\begin{array}}\newcommand{\ea}{\end{array}}
\newcommand{\beqa}{\begin{eqnarray}}\newcommand{\eeqa}{\end{eqnarray}}

\newtheorem{proposition}{Proposition}
\newtheorem{lemma}{Lemma}
\newtheorem{cor}{Corollary}

\begin{document}
\title{On rime Ansatz}
\author{ {Oleg Ogievetsky\footnote{Centre de Physique Th\'eorique, Luminy, 13288 Marseille, France} $\,$ and$\,$
Todor  Popov\footnote{Institute for Nuclear Research and Nuclear Energy, Sofia, 
1784, Bulgaria}}\\
} 
\date{}
\maketitle
\begin{abstract}\noindent
The ice Ansatz on matrix solutions of the Yang--Baxter equation is weakened to a condition which 
we call {\it rime}. Generic rime solutions of the Yang--Baxter equation are described. We prove 
that the rime non-unitary (respectively, unitary) R-matrix is equivalent to the Cremmer--Gervais 
(respectively, boundary Cremmer--Gervais) solution. Generic rime classical r-matices satisfy 
the (non-)homogeneous associative classical Yang-Baxter equation.
\end{abstract}

\noindent
Let $V$ be a vector space.
Among  solutions $\hat{R}\in\,$End$\, (V\otimes V)$ 
of the Yang-Baxter equation 
\beq
(\hat{R} \otimes 1 \!\! 1) (1 \!\! 1 \otimes \hat{R}) (\hat{R} \otimes 1 \!\! 1)=
(1 \!\! 1 \otimes \hat{R}) (\hat{R} \otimes 1 \!\! 1)(1 \!\! 1 \otimes \hat{R})\ ,\nonumber\eeq
\noindent there is a class characterized by the so called {\it ice condition} (see 
lectures \cite{QS} for details) which says that $\hat{R}^{i j }_{k l}$ can be different from zero
only if the sets of upper and lower indices coincide,
\[ \hat{R}^{i j }_{k l} \neq 0  \qquad \Rightarrow \qquad \{ i , j  \} \equiv \{ k , l \}.\]
\noindent We introduce the {\it rime} Ansatz, relaxing the ice condition: the entry 
$\hat{R}^{i j }_{k l}$ can be different from zero if the set of the lower indices is a subset of 
the set of the upper indices,
\[ \hat{R}^{i j }_{k l} \neq 0 \qquad \Rightarrow \qquad  \{ k,  l \} \subset \{ i , j  \}\ . \]
\noindent Matrices for which it holds will be referred to as {\it rime} matrices.
A rime matrix has the 
structure (to avoid redundancy, we fix $\beta_{ii}=0$, $\gamma_{ii}=0=\gamma'_{ii}$)
\beq\hat{R}_{kl}^{ij}=\alpha_{ij}\delta^i_l\delta^j_k 
+\beta_{ij}\delta^i_k\delta^j_l+\gamma_{ij} 
\delta^i_k\delta^i_l+\gamma'_{ij}\delta^j_k\delta^j_l\ \qquad\mbox{(no summation)}
\ .\label{rice}\eeq
\noindent We denote by $\alpha_i$ the 
diagonal elements $\hat{R}_{ii}^{ii}$, $\alpha_i:=\alpha_{i i }$. Throughout the text we shall 
assume that the matrix $\hat{R}$ is invertible which, in particular, implies that $\alpha_i\neq 0$ 
for all $i$. 

The ice and rime matrices are made of $4 \times 4$ elementary building blocks, respectively,
\beq\label{fri}\hat{R}^{ice}=
\left(\ba{cccc}\alpha_1&0&0&0\\0&\beta_{12}&\alpha_{12}&0\\
0&\alpha_{21}&\beta_{21}&0\\0&0&0&\alpha_2\ea\right)\qquad {\mathrm{and}}\qquad 
\hat{R}^{rime}=\left(\ba{cccc}\alpha_1&0&0&0\\ \gamma_{12}&\beta_{12}&\alpha_{12}&\gamma'_{12}\\
\gamma'_{21}&\alpha_{21}&\beta_{21}&\gamma_{21}\\0&0&0&\alpha_2\ea\right)\ .\eeq
\noindent We call a rime matrix ``strict" if  $\alpha_{ij}\gamma_{ij} \neq 0$ 
$\forall\ i$ and $j$, $i\neq j$. Strict rime matrices are necessarily not ice. 

All propositions hereafter are given without proofs, for further details see \cite{OP}.

\begin{proposition}$\!\!\!${\bf .} The Yang-Baxter equation implies that the quantity $\beta :=\beta_{ij}+\beta_{ji}$ does not depend
on $i$ and $j$. A strict  rime Yang-Baxter solution (\ref{rice})
has, up to a change of basis, the form
\beq\label{R}
\hat{R}_{k l }^{i j}(\beta_{ i j})=(1 - \beta_{j i} )\delta^i_l \delta^j_k + 
\beta_{i j} \delta^i_k \delta^j_l - 
\beta_{i j} \delta^i_k \delta^i_l + \beta_{ j i} \delta^j_k \delta^j_l \ . 
\eeq
\noindent
 \end{proposition} 


\begin{lemma}$\!\!\!${\bf .} \label{Heck}
The rime Yang-Baxter solution (\ref{R}) is of Hecke type, 
\beq\label{Hecke}\hat{R}^2(\beta_{ij})= \beta \hat{R}(\beta_{ij}) + (1-\beta)I \  ,\qquad I:=1\!\!1 \otimes 1\!\!1\ .\eeq
\noindent When $\beta\neq2$, $\hat{R}(\beta_{ij})$ is of $GL$-type: multiplicities 
of the eigenvalues
of $1$ and $\beta -1$ are, respectively, $\frac{n(n+1)}{2}$ and $\frac{n(n-1)}{2}$.
\end{lemma}

\begin{proposition}$\!\!\!${\bf .} \label{prop:rnu}
The non-unitary ($\beta\neq 0$) strict rime  Yang-Baxter  solutions  (\ref{R})  
are parameterized by a point $\phi\in \mathbb P \mathbb C^n$ 
in a projective space, $\phi=(\phi_1: \phi_2: \ldots :\phi_n)$, such that $\phi_i\neq 0$  and $\phi_i\neq \phi_j$, $i\neq j$. 
The $\beta_{ij}$'s of the solution $\hat{R}(\beta_{ij})$ are given by  
\beq \label{bphi}\beta_{ij}=\frac{\beta\phi_i}{\phi_i-\phi_j}\ , \qquad  
\hat{R}:= \hat{R}(\frac{\beta\phi_i}{\phi_i-\phi_j}) \ .
\eeq
\end{proposition}
\noindent If one and only one $\phi_i=0$, then $\hat{R}$ is again rime Yang-Baxter 
solution (it is not strict).

\begin{proposition} The unitary ($\beta =0$) rime R-matrices $\hat{R}(\beta^{(0)}_{ij})$ 
are parameterized by a vector 
$(\mu_1,\ldots ,\mu_n)$ such that $\mu_i \neq \mu_j$,
\beq\label{unitary}\beta_{ij}^{(0)} =\frac{1}{\mu_i - \mu_j}\ , \qquad \hat{R}_0:= \hat{R}(\beta^{(0)}_{ij}) \ .\eeq
\end{proposition}

The Cremmer-Gervais solution \cite{CG} of the Yang-Baxter equation in its general two-parametric 
form reads (we use a rescaled matrix with eigenvalues $1$ and $ -q^{-2}$)
\beq\label{CG}
(\hat{R}_{CG,p})^{ij}_{kl} = q^{-2\theta_{ij}} p^{i-j} \delta^{i}_{l} \delta^{j}_{k}
+ (1 - q^{-2}) \sum_{i \leq s < j} p^{i-s} \delta^{s}_{k} \delta^{i + j -s }_{l} 
- (1 - q^{-2}) \sum_{j < s < i} p^{i-s}   \delta^{s}_{k} \delta^{i + j -s }_l\ ,\eeq
\noindent where  $\theta_{ij}$ is the step function ($\theta_{ij} =1$ when $i >j$ and 
$\theta_{ij} =0$ when $i \leq j$). 

\begin{proposition} The rime R-matrix $\hat{R}$ (\ref{bphi}) transforms to the Cremmer-Gervais 
matrix
\beq\label{change}\hat{R}=(X(\phi)\otimes X(\phi))\,\hat{R}_{CG,1}\,\, (X^{-1}(\phi)\otimes X^{-1}(\phi)) \ , \qquad
\beta = 1 - q^{-2}\ ,\eeq
\noindent  by a change of basis  with the matrix $X(\phi)$ whose entries are the 
elementary symmetric polynomials  $e_i(x_1,\ldots ,x_N)=\mathop{\sum}\limits_{s_1< \ldots 
<s_{i}}x_{s_1}x_{s_2}\ldots x_{s_{i}}$,
\beq X^k_j (\phi) := e_{j-1} \, (\phi_{1},\ldots, \hat{\phi}_{k},
 \ldots,\phi_{n} ), \qquad 
\det X ={\prod}_{j<k} (\phi_{j}- \phi_{k})
 \ .
\label{matX}\eeq 
\end{proposition}

We now consider solutions of the classical Yang-Baxter (cYB) equation (i.e., classical r-matrices),
 $[r_{12},r_{23}] +[r_{12},r_{13}] +[r_{13},r_{23}]=0$ which are (quasi-)classical limits of 
the rime R-matrices above. Denote by $P$ the permutation operator,
 $P \, x\otimes y = y \otimes x$.

\begin{cor}The non-unitary rime R-matrix $R=P \hat{R}$ is linear in the parameter $\beta$,
\beq
\label{quant}R=I + \beta \, r\ ,  \quad \mbox{where} \quad
r= \sum_{i,j \, : \, i\neq j}  \frac{1}{\phi_i-\phi_j}(\phi_i e^i_j -\phi_j e^j_j)\otimes
(e^j_i - e^i_i) \ .
\eeq
\noindent $R$ is a
quantization of the non-skew-symmetric rime  r-matrix, $r + r_{21}=P -I$. Similarly 
$R_{CG,1}:= P\hat{R}_{CG,1} =I+ \beta r_{CG}$, thus the matrix $r$ is equivalent to the matrix 
\beq \label{clcr}r_{CG}=\sum_{i< j}\,\,
\sum_{s=1}^{j-i}(e^j_{i+s-1}\otimes e^i_{j-s+1}-e^i_{i+s-1}\otimes e^j_{j-s+1})\ ,\ r=
Ad_{X(\phi)}\otimes Ad_{X(\phi)} (r_{CG})\ .\eeq
\end{cor}

The non-skew-symmetric cYB solution such as $r$ and $r_{CG}$ are classified by 
Belavin-Drinfeld triples \cite{BD}. Gerstenhaber and Giaquinto introduced  the notion of  {\it 
boundary} skew-symmetric cYB solution: a solution which lies on the boundary of the space of  
skew-\-sym\-metric  solutions of the  modified classical Yang-Baxter equation \cite{GG} 
(these solutions  are in turn into one-to-one correspondence with non-skew-symmetric solutions of 
cYB). 

\begin{proposition} The unitary rime R-matrix ${R}_0=P\hat{R}_0$ is the 
quantization\footnote{The small parameter can be absorbed in the $\mu$'s} of the  
skew-sym\-met\-ric rime  matrix 
$r_0$
\beq
{R}_0 = 1\!\!\!1 + r_0 \ , \qquad r_0 = \sum_{i,j \, : \, i<j} \beta^{(0)}_{ij} Z^i_j \wedge Z^j_i 
\ ,
\eeq
\noindent where  $Z^i_j:= e^i_j - e^j_j$  generate a subalgebra of 
the matrix algebra and $x\wedge y:=x\otimes y-y\otimes x$.
The cYB solution $r_0$ is equivalent to the Cremmer-Gervais boundary solution \cite{GG,EH}
\beq
 b= \sum_{i<j} 
\sum_{k=1}^{j-i} e^{i+k}_{i} \wedge e^{j-k+1}_j\ , \qquad
r_0= Ad_{X(\mu)}\otimes Ad_{X(\mu)} (b) \ .
\eeq
\end{proposition}

By a sophisticated construction  \cite{GG} one can show  that $b$ is the boundary solution 
attached to $r_{CG}$ \cite{EH}. In the rime basis, the cYB solutions $b$ and $r_{CG}$ are 
transformed to the r-matrices $r_0$ and $r$, respectively. One of advantages of the rime basis is 
that the unitary limit is explicit: expand $\phi_i = 1 + \beta \mu_i + o(\beta) $ and take the 
limit
$$
\beta_{ij} = \frac{\beta(1 + \beta \mu_i+ o(\beta)) }{ \beta \mu_i - \beta \mu_j + o(\beta)}
\quad \stackrel{\beta \rightarrow 0}{\longrightarrow } \quad
\beta_{ij}^{(0)} =\frac{1}{\mu_i - \mu_j}\ .
$$
\noindent Hence ${R} \stackrel{ \beta \rightarrow 0}{\longrightarrow } {R}_0$ and  $\beta r \stackrel{ \beta \rightarrow 0}{\longrightarrow } r_0 $ (but the transform $X(\phi)$ becomes singular in the limit).

On the matrix level one can write $R_0$ as an exponential due to the nilpotency of $r_0$ \cite{EH}
$$r_0^2=0 \qquad \Rightarrow \qquad  {R}_0 = e^{ r_0} = I+ r_0 \ .
$$
\noindent The idempotency of $r$  yields a  similar exponential formula for $R$ 
$$
r^2 = - r  
\qquad \Rightarrow \qquad R= e^{h r} = I 
+ (1 - e^{-h})\, r\ ,
$$
\noindent 
so the quasi-classical approximation is exact in the renormalized  parameter 
$\beta=\beta(h)=1-e^{-h}$. 

Consider the splitting of the cYB equation: $ {\bf{A}}'(r) - {\bf A}(r) =0$, 
$${\bf A}(r):=r_{13}r_{12}-r_{12}r_{23}+r_{23}r_{13}\ , \qquad 
{\bf A}'(r):=r_{12}r_{13}-r_{23}r_{12}+r_{13}r_{23}\ .$$
\noindent  It was shown in \cite{MA1} that a solution of the associative  cYB (acYB) 
equation ${\bf A}({\bf r})=0$ defines 
the Newtonian coalgebra structure of an infinitesimal  
bialgebra \cite{JR}. The boundary r-matrices $b$ and $r_0$ satisfy  the acYB equations
$ {\bf A}(b)=0$, $ {\bf A}(r_0)=0$.

The non-skew-symmetric r-matrices $r_{CG}$ and $r$ satisfy the  non-homogeneous
acYB equation 
$${\bf A}(r)=-r_{13}\ .$$
\noindent 
and $r+r_{21}=P-I$ (thus ${\bf A}(\tilde{r})=\frac{1}{4} 1\!\! 1\otimes 1\!\!1 \otimes 1\!\!1 ,\quad\tilde{r}+\tilde{r}_{21}=P$, where $\tilde{r}=r+ \frac{1}{2}I$).

In the associative algebra generated by $r_{13}$, $r_{23}$, $r_{12}$ having as relations the 
non-homogeneous acYB equations ${\bf A}(r)=\xi r_{13}$, ${\bf A}'(r)=\xi r_{13}$ and the Hecke 
condition $r_{ij}^2=\xi r_{ij}+\eta$, $\{ ij\} =\{ 13\} ,\{ 23\} ,\{ 12\}$ ($\xi$ and $\eta$ are arbitrary constants)
the following identities hold \cite{OP}
$$
r_{12}r_{23}r_{12}= r_{23}r_{12}r_{23} \ ,\qquad\qquad r_{12}r_{13}r_{23}=r_{23}r_{13}r_{12}\ .
$$
\noindent i.e., the matrix $r$   satisfies the both forms of the ``quantum" Yang-Baxter equation.

The acYB (respectively, non-homogeneous acYB) solutions 
are related to the Rota--Baxter operators of zero (respectively, non-zero) weight. In the example
of integration and summation operators, the deformation from zero to non-zero weight is given
by the Euler-Maclaurin formula \cite{OS}. 
 
\noindent {\bf Acknowledgements.}
Todor Popov expresses his gratitude for the kind invitation and hospitality 
during the conference {\it Supersymmetry and Quantum Symmetries 07} in Dubna. The work 
was supported by the ANR project GIMP No.ANR-05-BLAN-0029-01.


\end{document}